\documentclass[11 pt]{amsart}
\usepackage{amsmath,enumerate,amsthm,amssymb,amscd}

\newcommand{\gothic}{\mathfrak}
\newcommand{\pd}{\operatorname{pd}}

\newcommand{\m}{{\gothic{m}}}

\newcommand{\Tor}{\operatorname{Tor{}}}

\renewcommand{\phi}{\varphi}

\newcommand{\lR}{\lambda_{R}}

\newtheorem{thm}{Theorem}

\newtheorem{prop}[thm]{Proposition}

\newtheorem{conjecture}[thm]{Conjecture}

\numberwithin{equation}{section}

\begin{document}

\title{A Note on a Conjecture of Watanabe and Yoshida}

\author {Lori McDonnell}

\address{Department of
Mathematics\\
University of Nebraska-Lincoln\\Lincoln,  NE 68588-0130}
\email{s-lmcdonn1@math.unl.edu}

\date{January 2010}
\
\bibliographystyle{amsplain}

\numberwithin{thm}{section}

\begin{abstract} We consider a conjecture of Watanabe and Yoshida in \cite{WY} concerning the Hilbert - Kunz multiplicity of an ideal in a Cohen-Macaulay ring and provide a proof of the conjecture in the case the ring is graded.
\end{abstract}

\maketitle

\section{Introduction}

The purpose of this note is to examine a conjecture of Watanabe and Yoshida and to provide a proof of  this conjecture in the case we have a graded ring $R$.  We begin with some notation.  Throughout this note all rings are assumed to be commutative Noetherian of prime characteristic $p>0$ unless otherwise noted.  Let $(R,m)$ be a local ring of dimension $d$ with maximal ideal $\mathfrak{m}$.  For an ideal $I$ of $R$ and $q=p^{e}$ (for some $e$), we let $I^{[q]}$ denote the ideal of $R$ generated by the set $\{i^{q} : i \in I \}$.
Given an $m$-primary ideal $I$ of $R$, let $\lambda_{R}(R/I)$ denote the length of $R/I$ as an $R$-module.  One then defines the Hilbert - Kunz multiplicity of $I$ by
$$
e_{HK}(I,R) := \lim_{q \rightarrow \infty} \frac {\lambda_{R}(R/I^{[q]})} {q^{d}}.
$$
It was shown by Monsky \cite{MO} that this limit exists and is positive for all such ideals.

In \cite{WY}, Watanabe and Yoshida make the following conjecture:

\begin{conjecture} \label{main conj}
Let $\left( R,m\right) $ be a Cohen-Macaulay local ring of characteristic $p>0$. \ Then
\begin{enumerate}
\item For any $m$-primary ideal $I$, one has $e_{HK}(I, R) \geq
\lambda_{R}(R/I)$.
\item For any $m$-primary ideal $I$ with $\pd_{R}\left( R/I\right)
\,<\infty$, one has $e_{HK}\left( I,R\right) =\lambda _{R}\left(
R/I\right)$. \end{enumerate}
\end{conjecture}

Dutta \cite{DU} has shown Conjecture \ref{main conj} holds when $R$ is a complete
intersection ring (but not necessarily graded).  However, neither part of the conjecture holds for all Cohen-Macaulay rings.  Miller and Singh \cite{MS} have constructed a module $M$ of finite projective dimension over a Gorenstein ring $R$ of dimension $5$ with \[220 = \lim_{n \rightarrow \infty}\frac{\lR(F^{n}_{R}(M))}{p^{5n}} < \lR(M)=222.\]  From the proof of Theorem 6.4 in \cite{KK}, there exist  $\m$-primary ideals $J, I_{1}, \dots, I_{t}$ of finite projective dimension such that $I_{1}, \dots, I_{t}$ are parameter ideals and 
\begin{enumerate}
\item $\lR(M)=\lR(R/J)-\sum_{i=1}^{t}\lR(R/I_{i})$, and
\item $\lim_{n\rightarrow \infty}\frac{\lR(F^{n}(M))}{p^{5n}}=e_{HK}(J, R)-\sum_{i=1}^{t}\lR(R/I_{i})$.
\end{enumerate}
Thus, for the ideal $J \subset R$, we have $e_{HK}(J, R)<\lR(R/J)$, contradicting both parts of the conjecture above.

In this note all graded rings $R=\oplus_{i\geq0}R_{i}$ are finitely generated over the Artinian local ring $R_{0}$.  What we are able to prove, using basic properties of Poincar\'e series and a result of Avramov and Buchweitz \cite{AB}, is the following:

\begin{thm} \label{main thm}
Let $R$ be a graded ring of characteristic $p>0$ and dimension $d$, and $I$ a homogeneous ideal with $\lambda(R/I)<\infty$ and $\pd(R/I)<\infty$.  Then for every $q=p^{e}$, one has $\lambda(R/I^{[q]})=q^{d}\lambda(R/I)$.  In particular, $e_{HK}(I, R)=\lambda(R/I)$.
\end{thm}

Theorem \ref{main thm} also follows from a conjecture of Szpiro \cite[Conjecture C2]{S}. Szpiro sketches a proof of this conjecture in the graded case, using different methods than what we employ here.

\section{Proof of Theorem \ref{main thm}}

Let $R$ be a graded ring and $M$ a nonzero finitely generated graded $R$-module. Let $P_{M}\left( t\right) =\sum_{i\in \mathbb{Z}}\lambda_{R_{0}} \left( M_{i}\right) t^{i}$ denote the Hilbert series for $M$.
Note that if $\lambda_{R} \left( M\right) <\infty ,$ we have $\lambda_{R} \left(
M\right) =P_{M}\left( 1\right) .$ Further note $P_{M\left( -k\right)
}\left( t\right) =t^{k}P_{M}\left( t\right) $ for any $k\in \mathbb{Z}$, where if $M= \oplus_{i\in \mathbb{Z}}M_{i}$, $M(-k)$ denotes the graded $R$-module with $M(-k)_{i}=M_{i-k}$.   We recall the following proposition concerning the Hilbert series for $M$.

\begin{prop} \label{thm: P as fraction}\textnormal{(cf.} \cite[Theorem 5.5]{SM} \textnormal{)}
Let $R$ be a graded Noetherian ring and $M$ a nonzero finitely generated graded $R$-module of dimension $l$.  Then there exist positive integers $s_{1,}\dots ,s_{l}$ and a
polynomial $p\left( t\right) \in \mathbb{Z}
\left[ t,t^{-1}\right] $ with $p\left( 1\right) \neq 0$ such that
\begin{equation*}
P_{M}\left( t\right) =\frac{p\left( t\right) }{\prod_{i=1}^{l}\left(
1-t^{s_{i}}\right) }.
\end{equation*}
\end{prop}

The following result is a special case of a result in \cite{AB}.  Since the proof is not difficult, we include it here for completeness.

\begin{prop}  \cite[Lemma 7]{AB} \label {lem: P and chi}
Let $M$ be a finitely generated graded $R$-module with finite length and
finite projective dimension. \ Let $\mathfrak{m}$ be the maximal ideal of $R_{0}$ and
let $K=R/(\mathfrak{m}+R_{+})=R_{0}/\mathfrak{m}$. Set $\chi _{M}^{R}\left( t\right)
:=\sum_{i}\left( -1\right) ^{i}P_{\Tor_{i}^{R}(M,K)}\left( t\right) .$ \ Then
one has $P_{M}\left( t\right) =\chi _{M}^{R}\left( t\right) P_{R}\left(
t\right)$.

\begin{proof}
Let $F_{\cdot }:0\rightarrow F_{s}\rightarrow
\cdots \rightarrow F_{1}\rightarrow F_{0}\rightarrow 0$ \ be
a minimal free resolution for $M$ where $F_{i}\cong \oplus _{j=
0}^{r_{i}}R\left( -j\right) ^{b_{ij}}$ with $b_{ij}\in \mathbb{N}
$. \ Then evaluating the Hilbert series for $M$ using the resolution,
one gets
\begin{equation}
P_{M}\left( t\right) =\sum_{i,j\in \mathbb{Z}
}\left( -1\right) ^{i}b_{ij}P_{R}\left( t\right) t^{j}.  \label{MtoR}
\end{equation}
Note this sum is well-defined as there are only finitely many nonzero $b_{ij}$.  Moreover, $\chi _{M}^{R}\left( t\right) $ is also well-defined as $\Tor^{R}_{i}(M, K)$ is finitely generated graded and $\Tor^{R}_{i}(M, K)=0$ for $i>\pd_{R}M$.
Now, tensoring $F_{\cdot }$ with $K$, we obtain the complex
\[\cdots \rightarrow \oplus K\left( -j\right) ^{b_{sj}}\rightarrow \cdots \rightarrow
\oplus K\left( -j\right) ^{b_{0j}}\rightarrow  0\]
where each map is the zero map.
Note that for $n \in \mathbb{Z}$, $\Tor_{i}^{R}\left( M,K\right) _{n}=H_{i}\left( F_{\cdot }\otimes
K\right) _{n},$ so $\dim _{K}\Tor_{i}^{R}\left( M,K\right) _{n}= \dim
_{K}\left( \left( \oplus _{j\in \mathbb{Z}
}K\left( -j\right) ^{b_{ij}}\right) _{n}\right) =\sum_{j\in \mathbb{Z}
}b_{ij}\dim _{K}K_{n-j}=b_{in}<\infty $ \ (as $K_{i}=0$ for $i\neq 0$).  Thus, $P_{\Tor_{i}^{R}(M,K)}(t)=\sum_{j\in \mathbb{Z}}b_{ij}t^{j}$.

Now, by the invariance of Euler-Poincare characteristics, we have
\begin{eqnarray*}
\chi _{M}^{R}\left( t\right) &=& \sum_{i\in \mathbb{Z}}\left( -1\right) ^{i}P_{\Tor_{i}^{R}\left( M,K\right) }\left( t\right)  \\
&=&\sum_{i,j\in \mathbb{Z}
}\left( -1\right) ^{i}b_{ij} t^{j} \\
&=&\frac{P_{M}\left( t\right)} {P_{R}\left(
t\right) }  .\text{ \ \ \ \ \ \ \ \ \ \ \ \ \ (by equation (\ref{MtoR}))} \\
\end{eqnarray*}
This gives
$P_{R}\left( t\right) \chi _{M}^{R}\left( t\right) =P_{M}(t)$.  Note also that we have $\chi^{R}_{M}(t) \in \mathbb{Z}[t,t^{-1}]$.
\end{proof}
\end{prop}

Note that if we let $F$ denote the Frobenius functor, then the Hilbert-Kunz multiplicity of an ideal is found by taking a limit of increasing iterations of the Frobenius of the $R$-module $R/I$, i.e. $e_{HK}(I, R)= \lim_{e \rightarrow \infty} F^{e}(R/I)/p^{ed}$.  Now, if we can show $\lambda_{R} (F^{e}(R/I)) =p^{ed} \lambda_{R} (R/I)$ for all $e$ sufficiently large, then  $e_{HK}(I,R) = \lambda_{R} (R/I)$.
We now prove the main theorem.

\begin{thm}
Let $R$ be a graded ring of characteristic $p$ and let $M$ be a finitely generated graded $R$-module with $\lambda_{R}(M)<\infty $ and $\pd_{R}(M)<\infty$.  Then $\lambda_{R} \left( F^{e}(M) \right) =q^{d}\lambda_{R} \left(M\right)$ for all $q=p^{e}$.  In particular,
one has $e_{HK}(I, R) =\lambda_{R} \left( R/I\right)$ for all zero-dimensional homogeneous ideals $I$ of finite projective dimension.

\begin{proof}
\ Let $G_{\cdot }$ be a minimal graded free resolution for $M,$ with $G_{i}=\overset
{r_{i}}{\underset{j=0}{\oplus }}R\left( -j\right) ^{b_{ij}}$ and $
b_{ij}\in \mathbb{N}$.  Let $F^{e}\left( -\right) $ denote the Frobenius functor and $q=p^{e}.$
\ Note that $L_{\cdot }=F^{e}\left( G_{\cdot }\right) $ is a minimal graded free
resolution for $F^{e}(M)$ by \cite[Theorem 1.7]{PS2} where each twist by $j$ in $G_{\cdot}$ is multiplied by a factor of $q$ and the $b_{ij}$ remain the same. That is, $L_{i}=\overset{r_{i}}{\underset{j=1}{\oplus}}R\left( -jq\right) ^{b_{ij}}.$ \ Now, by the lemma above, we
have $P_{M}\left( t\right) =\chi _{M}^{R}\left( t\right) P_{R}\left(
t\right) $ and $P_{F^{e}(M)}\left( t\right) =\chi _{F^{e}(M)}^{R}\left( t\right) P_{R}\left( t\right) .$ \

By Proposition \ref{thm: P as fraction}, we can write $P_{R}\left( t\right) =\dfrac{p\left(
t\right) }{\prod_{i=1}^{d}\left( 1-t^{s_{i}}\right) }$ where $d=\dim R$ and the
$s_{i}\in \mathbb{N}$ are nonzero. Each term in the denominator can be factored as $1-t^{s_{i}}=(1-t) g_{i}\left( t\right) $ with $g_{i}(t) \in \mathbb{Z} \left[ t\right] $ and $g_{i}(1) \neq 0$.  Letting $g(t) = \prod_{i} g_{i}(t)$,
we can rewrite $P_{R}\left( t\right) =\dfrac{p\left( t\right) }{\left(
1-t\right) ^{d}g\left( t\right) }$ where $p(1)/g(1)$ is a nonzero rational number.

Since $\lambda_{R} \left(M\right) <\infty ,$ $P_{M}\left( t\right) \in \mathbb{Z}
\left[ t, t^{-1}\right] $.  So, we have
\begin{equation*}
P_{M}\left( t\right) =\chi _{M}^{R}\left( t\right) P_{R}\left( t\right)
=\chi _{M}^{R}\left( t\right) \dfrac{p\left( t\right) }{\left( 1-t\right)
^{d}g\left( t\right) }\in \mathbb{Z}[ t, t^{-1}] .
\end{equation*}
\ Since $p\left( 1\right) \neq 0,$ we must have $\left( 1-t\right) ^{d}$ divides
$\chi _{M}^{R}\left( t\right)$ in $\mathbb{Z}[t, t^{-1}]$; say $\chi _{M}^{R}\left( t\right) =
\tilde{\chi}_{M}^{R}\left( t\right) \cdot \left( 1-t\right) ^{d}$, for some $\tilde{\chi}_{M}^{R}(t) \in \mathbb{Z}[t, t^{-1}]$.  So,
\begin{equation}
P_{M}\left( t\right) =\tilde{\chi}_{M}^{R}\left( t\right) \frac{p\left(
t\right) }{g(t) }.  \label{PM}
\end{equation}

In the proof of Proposition \ref{lem: P and chi}, we see that $\chi _{M}^{R}\left( t\right)
=\sum_{i,j\in \mathbb{Z}}\left( -1\right) ^{i}b_{ij}t^{j}$.  Applying this to the resolution $L_{\cdot}$ of $F^{e}(M)$, we also have $\chi _{F^{e}(M)}^{R}\left( t\right)
=\sum_{i,j\in \mathbb{Z}}\left( -1\right) ^{i}b_{ij}t^{qj}=\chi _{M}^{R}\left( t^{q}\right) =
\tilde{\chi}_{M}^{R}\left( t^{q}\right) \cdot \left( 1-t^{q}\right) ^{d}.$
\ Thus,
\begin{eqnarray*}
P_{F^{e}(M)}\left( t\right) &=&\chi _{F^{e}(M)}^{R}\left( t\right) \frac{p\left( t\right) }{\left( 1-t\right) ^{d}g\left(
t\right) } \\
&=&\tilde{\chi}_{M}^{R}\left( t^{q}\right) \left( 1-t^{q}\right) ^{d}\frac{
p\left( t\right) }{\left( 1-t\right) ^{d}g\left( t\right) } \\
&=&\tilde{\chi}_{M}^{R}\left( t^{q}\right) \left( 1+t+\cdots
+t^{q-1}\right) ^{d}\frac{p\left( t\right) }{g\left( t\right) }.
\end{eqnarray*}
Letting $t=1$, we get
\begin{eqnarray*}
\lambda_{R} (F^{e}(M)) &=&P_{F^{e}(M)}\left( 1\right) =\tilde{\chi}_{M}^{R}\left( 1^{q}\right) \left(
1+1+\cdots +1^{q-1}\right) ^{d}\frac{p\left( 1\right) }{g\left( 1\right) } \\
&=&\tilde{\chi}_{M}^{R}\left( 1\right) \left( q\right) ^{d}\frac{p(1)}{%
g\left( 1\right) } \\
&=&q^{d}P_{M}\left( 1\right) \text{ \ \ \ \ \ \ (by (\ref{PM}))} \\
&=&q^{d}\lambda_{R} \left( M\right).
\end{eqnarray*}
Finally, if $I$ is a zero-dimensional homogeneous ideal and $q=p^{e}$, we have $e_{HK}(I, R) =\lim_{q\rightarrow \infty }\lambda_{R}
\left(F^{e}(R/I)\right) /q^{d}=\lambda_{R} \left( R/I\right) .$
\end{proof}
\end{thm}

%In light of this result, one might ask the stronger question related to Conjecture \ref{main conj}:

%\begin{ques} \label{ques}
%Let $(R,\mathfrak{m})$ be Cohen-Macaulay and $M$ a finitely generated $R$-module such that $\lambda_{R}(M)<\infty$ and $\pd_{R}(M)<\infty$.  Under what conditions is $\lambda_{R}(F^{n}(M))=p^{nd}\lambda_{R}(M)$ for all $n$?
%\end{ques}

%Dutta has shown that Question \ref{ques} has an affirmative answer when $R$ is a complete intersection; in fact, Miller \cite{MI} has extended Dutta's result to show that $\lambda_{R}(F^{n}(M))=p^{nd}\lambda_{R}(M)$ for all $n$ if and only if $\pd_{R}(M)<\infty$.  She also has shown that $\lambda_{R}(F^{n}(M))=p^{nd}\lambda_{R}(M)$ if and only if $\pd_{R}(M)<\infty$ when $R$ is a Cohen-Macaulay local ring of dimension one which is geometrically unibranch.  Dutta has further shown that Question \ref{ques} holds true when $R$ is Gorenstein of dimension less than or equal to three. However, Miller and Singh \cite{MS} have constructed a module over a Gorenstein ring of dimension five to show that the equality $\lambda_{R} \left(F^{n}(M\right) )=p^{nd}\lambda_{R} \left( M\right) $ for all $n$ does not hold in general over Gorenstein rings.

\section*{Acknowledgments}
The author would like to thank her advisor Tom Marley for his direction and guidance in this research and the editing of this note.  Thanks also to Luchezar Avramov for directing us to \cite{AB} and the methods used there.

\end{document}